\documentclass[thmsa]{article}
\usepackage{amssymb}

\usepackage{sw20lart}

\input{tcilatex}
\begin{document}

\author{\textit{Stancho Dimiev}}
\title{Pseudogroup structures on Spencer manifolds}
\date{\textit{e-mail address: sdimiev@math.bas.bg}}
\maketitle

\begin{abstract}
This note is a continuation of the paper [2] (see references). We describe
some natural pseudogroup structures on almost complex manifolds of type $m$.
A kind of coherency is discussed for the sheaf of almost holomorphic
functions.
\end{abstract}

\section{Recall of definitions}

Let $(M,J)$ be an almost complex manifold equipped with a smooth coordinates 
$(x^j)$, $j=1$, $2$, $....$, $2n$. In [1] complex self-conjugate coordinates 
$(z^j$, $\overline{z}^j)$, $j=1$, $2$, $....$, $n$, are used, where $%
z^j=x^{2j-1}+ix^{2j}$ and $\overline{z}^j=x^{2j-1}-ix^{2j}$. A function $%
f:U\rightarrow \mathbf{C}$, where $U$ is an open subset of $M$, is called
almost holomorphic if $J^{*}df=idf$. Here by $J^{*}$ the action of $J$ on
differential forms of $M$ is denoted, i.e. $J^{*}(\omega )X:=\omega (X)$,
where $\omega $ is a differential form on $M$, and $X$ is a vector field on $%
M$.

Let $m$ be an integer, $m\leq n$. We say that a local Spencer coordinate
system is defined on $(M,J)$, iff

(1) there exist an open subset $U$ of $M$ and $m$ different functionally
independent almost holomorphic functions $f_j:U\rightarrow \mathbf{C}$, $j=1$%
,$...$, $m$, such that

(2) the sequence $f_1$, $...$, $f_m$ is a maximal sequence of functionally
independent on $U$ almost-holomorphic functions

(3) the sequence $(U,\,\,f_1$, $...$, $f_m$, $z^{m+1}$, $...$, $z^n$, $%
\overline{f}_1$, $...$, $\overline{f}_m$, $\overline{z}^{m+1}$, $...$, $%
\overline{z}^n)$ determines a local self-conjugate coordinate system on $%
(M,J)$.

The condition (3) can be reformulate in terms of smooth coordinates as
follows:

(3')The sequence $(U,\,\,\func{Re}f_1$, $\func{Im}f_1$, $...$,\thinspace $%
\func{Re}$ $f_m$, $\func{Im}f_m$, $x^{m+1}$, $y^{m+1}$, $...$, $x^{2n}$, $%
y^{2n})$ is a local coordinate system in smooth coordinates on $(M,J)$.

Setting $u_j=\func{Re}f_j$, $v_j=\func{Im}f_j$, we have $dv_j=J^{*}du_j$, $%
j=1$, $...$, $m$. These are Cauchy-Riemann equations for $f_j$.

Clearly, each Spencer coordinate system on $M$ defines a diffeomorphism of $%
U $ in $\mathbf{C}^n\times \overline{\mathbf{C}}^n$, or equivalently in $%
\mathbf{R}^{2n}$. Complexifying the first $2m$ coordinate we can take $%
\mathbf{C}^m\mathbf{\times R}^{2n-2m}$, i.e. $\mathbf{R}^{2n}\equiv \mathbf{C%
}^m\times \mathbf{R}^{2n-2m}$.

We shall consider the smooth submersion $\mathbf{R}^{2n}\rightarrow \mathbf{C%
}^m$ defined as a composition of the above mentioned diffeomorphism of $U$
on $\mathbf{R}^{2n}$ and the projection of $\mathbf{R}^{2n}$ on $\mathbf{C}%
^m $%
\[
(f_1,...,f_m,z^{m+1},...,z^n,\overline{f}_1,...,\overline{f}_m,\overline{z}%
^{m+1},...,\overline{z}^n)\rightarrow (f_1,...,f_m)\text{ .} 
\]

This composition is denoted by $f_U$ [more precisely by $(f_1,...,f_m)_U$]
and the image of $U$ by $f_U$ is denoted by $U^c\,_m$.

As it was proved in [2], each almost-holomorphic function $h:U\rightarrow 
\mathbf{C}$ is a superposition of a holomorphic function $%
H:U^c\,_m\rightarrow \mathbf{C}$ and the almost-holomorphic functions $f_j$, 
$j=1$, $...$, $m$, i.e. $h=H(f_1,...,f_m)$, where $H=H(w_1,...,\,w_m)\in 
\mathbf{O}(U^c\,_m)$.

It was proved also that for each two systems of almost-holomorphic functions 
$(f_1,...,f_m)$ and $(g_1$, $...$, $g_m)$ defined on $U\cap V$, there exists
a bijective holomorphic transition mapping between the corresponding Spencer
coordinate systems.

As a consequence one obtain a theorem formulated by spencer [1], namely

\textbf{Theorem}. On each paracompact almost-complex manifold $(M,J)$ of
constant Spencer type $m$ there exists a locally finite covering $\{U_\alpha
\}$ constituted by self-conjugate Spencer's coordinate systems $(U_\alpha $, 
$w^j\,_\alpha $, $...)$, $j=1$, $...$, $m$, such that on every intersection $%
U_\alpha \cap V_\beta $ the Spencer coordinates $(w^j\,_\alpha )$ change
biholomorphically in the other Spencer's coordinates $(w^j\,_\beta )$.

\section{Pseudogroup structures}

By $X$ a topological space is denoted, respectively differentiable manifold $%
(X=M)$, or vector space $(X=\mathbf{R}^n$, $X=\mathbf{C}^n)$. We recall the
notion of pseudogroup $\Gamma $ of local homeomorphisms of $X$, respectively
- local diffeomorphisms of $M$, or local biholomorphisms of $\mathbf{C}^m$.
Let $U$ e an open subset of $X$, respectively of $M$, or of $\mathbf{C}^m$,
and $\varphi :U\rightarrow X$ be a local mapping (homeomorphism,
diffeomorphism, biholomorphism). The open set $U$ is called the source of $%
\varphi $, and the image of $\varphi (U)$ (which is also an open set in $X$)
is called a target of $\varphi $.

A family of local mappings $\{(U,\varphi ):$ $U$ varying in a part of the
set of open subsets of $X\}$ is by definition a pseudogroup of local
mappings (homeomorphisms, diffeomorphisms, biholomorphisms) if the following
axioms are valid:

(1) If $(U,\varphi )$, $(V,\psi )\in \Gamma $ and $\varphi (V)\subset U$,
then $(V,\varphi \circ \psi )\in \Gamma $,

(2) If $(U,\varphi )\in \Gamma $, and $V=\varphi (U)$, then $(V,\,\varphi
^{-1})\in \Gamma $,

(3) If $(U,\varphi )\in \Gamma $, and $V\subset U$, then $(V,\,\varphi \mid
V)$ (the restriction remains in $\Gamma $),

(4) If $(U,\varphi )\in \Gamma $ and every point of $U$ admits a
neighborhood on which the restriction of $\varphi $ is in $\Gamma $, then $%
\varphi $ is in $\Gamma $,

(5) The restriction of the identity on every source of an element of $\Gamma 
$ is in $\Gamma $.

Let $(X,\Gamma _X)$ be a topological space $X$ equipped with a pseudogroup
of local homeomorphisms $\Gamma _X$, and $(Y,\Gamma _Y)$ is another
topological space equipped with a pseudogroup of local homeomorphisms $%
\Gamma _Y$. We say that $(X,\,\Gamma _X)$ is defined over $(Y,\,\Gamma _Y)$
if for every source $U$ of an element of $\Gamma _X$, $(U,\,\varphi )\in
\Gamma _X$, there exists continuous mappings $f_U$, $f_{\varphi (X)}$ and $%
\psi $such that the following diagram should be commutative

\begin{center}
$
\begin{array}{lllll}
U &  & \stackrel{\varphi }{\rightarrow } &  & \varphi (U) \\ 
&  &  &  &  \\ 
\downarrow f_U &  &  &  & \downarrow f_{\varphi (U)} \\ 
&  &  &  &  \\ 
V &  & \stackrel{\psi }{\rightarrow } &  & \psi (V)
\end{array}
$,
\end{center}

where $V=f_U(U)$ is a source of $(V,\,\psi )\in \Gamma _Y$ and $\psi
(V)=f_{\varphi (U)}(\varphi (U))$, i.e. 
\[
f_{\varphi (U)}\circ \varphi =\psi \circ f_U\text{ \thinspace \thinspace
\thinspace .} 
\]

By $\Gamma _d(M)$ the transitive pseudogroup of all local diffeomorphisms of
the differentiable manifold $M$ is denoted. In the case $M$ is an
almost-complex manifold, $(M,\,J)$, we shall consider local
almost-holomorphic diffeomorphisms of $M$, i.e. the diffeomorphisms $%
f:U\rightarrow M$, \thinspace $U$ being an open subset of $M$, and $f$
satisfying the condition 
\[
f_{*}\circ J=J\circ f_{*}\text{ ,} 
\]

where $f_{*}$ is the tangent mapping (the differential) of the mapping $f$.

The minimal transitive pseudogroup of local almost-holomorphic mappings will
be denoted $\Gamma _{ahd}(M)$. This pseudogroup is a subpseudogroup of $%
\Gamma _d(M)$. It is to recall here that the composition of two
almost-holomorphic mappings is also an almost-holomorphic mapping where it
is defined. The same is true for the inverse of an almost-holomorphic
mapping.

Accordingly to the theorem formulated in previous paragraph, we shall
consider all local almost-holomorphic diffeomorphisms $\Phi _{\alpha \beta
}:U_\alpha \rightarrow U_\beta $, and corresponding biholomorphisms $\varphi
_{\alpha \beta }:(U_\alpha )^c\,_m\rightarrow (U_\beta )^c\,_m$ in $\mathbf{C%
}^m$, where $\{U_\alpha \}$ are the sources of local Spencer coordinates
systems. The family $\{\Phi _{\alpha \beta }\}$ generates a subpseudogroup
of $\Gamma _{ahd}(M)$ which will be denoted by $\Gamma _{spd}(M)$. Denoting
by $\Gamma _h(\mathbf{C}^m)$ the pseudogroup of all local biholomorphisms in 
$\mathbf{C}^m$, we consider the family $\{\varphi _{\alpha \beta }\}$ and
the generated subpseudogroup of $\Gamma _h(\mathbf{C}^m)$ which will be
denoted by $\Gamma _{sph}(\mathbf{C}^m)$. The pseudogroup $\Gamma _{spd}(M)$
is over the pseudogroup $\Gamma _{sph}(\mathbf{C}^m)$ according to the above
introduced definition. The pseudogroups $\Gamma _{spd}(M)$ and $\Gamma
_{sph}(\mathbf{C}^m)$ are called Spencer pseudogroup of the almost complex
manifold $(M,\,J)$ of type $m$. This means that for every $\Phi
:U\rightarrow V$, $(U,\,V)\in \Gamma _{spd}(M)$, we have $f_V\circ \Phi
=\varphi \circ f_U$, where $(U^c\,_m,$ $\varphi )\in \Gamma _{sph}(\mathbf{C}%
^m)$ and also $\Phi _{*}\circ J=J\circ \Phi _{*}$. The following diagram is
commutative

\begin{center}
$
\begin{array}{lllll}
U &  & \stackrel{\Phi }{\rightarrow } &  & V \\ 
&  &  &  &  \\ 
\downarrow f_U &  &  &  & \downarrow f_V \\ 
&  &  &  &  \\ 
U^c\,_m &  & \stackrel{\varphi }{\rightarrow } &  & V^c\,_m
\end{array}
\,$.
\end{center}

\section{$\Gamma -$ manifolds, integrability of $(M$, $J)$}

Let $\Gamma $ be a pseudogroup of local diffeomorphisms of one
differentiable manifold $M$. It is said that a $\Gamma $-structure on the
manifold $M$ is defined if an atlas of local coordinates systems is
introduced in such a way that all transition transformations between them to
belong to the pseudogroup $\Gamma $.

Let $(M,\,J)$ be an almost-complex manifold of type $m$. We shall consider
the corresponding Spencer pseudogroups $\Gamma _{spd}(M)$ and $\Gamma _{sph}(%
\mathbf{C}^m)$. Let us suppose that $N$ is an orientable differentiable $%
(2m) $-manifold equipped with a $\Gamma _{sph}(\mathbf{C}^m)$-structure. We
remark that the problem of the existence of $\Gamma $-structures on a given
manifold, especially $\Gamma _{sph}(\mathbf{C}^m)$-structures, is a
difficult problem.

Denoting $f=\{f_U\}$ as the family of all local $m$-projections $%
f_U:U\rightarrow \mathbf{C}^m$ defined by almost-holomorphic coordinates $%
f_1 $, $f_2$, $...$, $f_m$, we obtain the following diagram

\begin{center}
$
\begin{array}{lllll}
M &  & \equiv &  & M \\ 
&  &  &  &  \\ 
\downarrow f &  &  &  & \downarrow F \\ 
&  &  &  &  \\ 
\mathbf{C}^m &  & \stackunder{atlas}{\stackrel{\theta }{\leftarrow }} &  & N
\end{array}
\,$,
\end{center}

where $F$ is defined locally as follows: 
\[
F(U)=\theta ^{-1}\circ f(U)\text{ ,} 
\]

$U$ being an open subset of $M$.

In the case $m=n$ [$M$ is a $(2n)$-manifold] the $\Gamma _{sph}(\mathbf{C}%
^n) $-structure on $N$ is a structure of a complex manifold, and $F$ is a
diffeomorphism. So the almost complex manifold $(M,\,J)$ is diffeomorphic to
a complex manifold $N$. Taking $N=\mathbf{C}^n$ as differentiable manifold,
we obtain that $\Gamma _{sph}(\mathbf{C}^n)$ defines an atlas of
biholomorphisms on $M$. This implies that $(M,\,J)$ is an integrable
manifold.

In the case of $4$-dimensional almost-complex manifold $M$ of type $m=1$
there are Spencer pseudogroup structures $\Gamma _{spd}(M)$ and $\Gamma
_{sph}(\mathbf{C})$ The corresponding $2$-dimensional manifold $N$, equipped
with a $\Gamma _{sph}(\mathbf{C})$ must be a Riemann surface.

\section{$\Gamma $- coherency}

The notion of coherent sheaf is local. Having a pseudogroup $\Gamma $ on a
manifold $M$ we can assign to each source $U$ of $\Gamma $ the set $\Gamma
(U)$ of different $\Gamma $-objects (functions, vector fields) defined on $U$%
. The mapping $U\rightarrow \Gamma (U)$ defines a sheaf.

Our purpose is to discuss a kind of coherency on almost complex manifolds
using the introduced Spencer pseudogroup structures of type $m$. It is not
difficult to see that almost-holomorphic functions on $(M,\,J)$ define a
sheaf. We denote this sheaf by $O_{ah}(M)$ and, respectively, by $O_h(%
\mathbf{C}^m)$ - the sheaf of all holomorphic functions on $\mathbf{C}^m$.
According to the famous Oka theorem $O_h(\mathbf{C}^m)$ is a coherent sheaf,
but the same is not true in general for the subsheaf $O_h(U)$, where $U$ is
an open subset of $\mathbf{C}^m$.

\textbf{Proposition}. For every source $U\in \Gamma _{spd}(M)$, the sheaf $%
O_{ah}(U)$ is a coherent sheaf of almost holomorphic functions if and only
if the corresponding sheaf of holomorphic functions $O_h(U^c\,_m)$ is
finitely generated as a subsheaf of $O_h(\mathbf{C}^m)$.

Proof. It is enough to remark that $O_{ah}(U)$ is an inverse image of $%
O_h(U^c\,_m)$. Being a finitely generated subsheaf of $O_h(\mathbf{C}^m)$
the sheaf $O_h(U^c\,_m)$ is a coherent sheaf too. $\square $

\section{Perspectives}

It seems that one can develope a deformation theory for almost complex
manifolds of type $m$ following some ideas of Donald Spencer (see [2]).

\end{document}